\documentclass[12pt, oneside,reqno]{amsart}   	
\usepackage{geometry}                		
\usepackage{graphicx}				
\usepackage{amssymb}
\usepackage{amsmath}
\usepackage[all]{xy}
\usepackage{mathtools}
\usepackage{tikz-cd}
\usepackage{enumitem}

\usepackage{tikz-cd}

\newtheorem{thm}{Theorem}[section]
\newtheorem{lem}[thm]{Lemma}
\newtheorem{cor}[thm]{Corollary}
\newtheorem{prop}[thm]{Proposition}

\newtheorem{rem}{Remark}[section]

\newtheorem{defn}{Definition}[section]

\numberwithin{equation}{section}

\newcommand{\Balpha}{\mbox{$\hspace{0.12em}\shortmid\hspace{-0.62em}\alpha$}} 

\def\Pb{\ifmmode{\Bbb P}\else{$\Bbb P$}\fi}
\def\Z{\ifmmode{\Bbb Z}\else{$\Bbb Z$}\fi}
\def\C{\ifmmode{\Bbb C}\else{$\Bbb C$}\fi}
\def\R{\ifmmode{\Bbb R}\else{$\Bbb R$}\fi}
\def\S{\ifmmode{S^2}\else{$S^2$}\fi}

\def\S{\cal S}

\newenvironment{pf}{\paragraph{Proof:}}{\hfill$\square$ \newline}

\begin {document}
	
\title[Regularity and Stability Results for LSF via MCF with Surgery]{Regularity and Stability Results for the Level Set Flow via the Mean Curvature Flow with Surgery}
\begin{abstract}In this article we use the mean curvature flow with surgery to derive regularity estimates for the level set flow going past Brakke regularity in certain special conditions allowing for 2-convex regions of high density. We also show a stability result for the plane under the level set flow.  \end{abstract}

\author {Alexander Mramor}
\address{Department of Mathematics, University of California Irvine, CA 92617}
\email{mramora@uci.edu}

\maketitle

\section{Introduction.}

The mean curvature flow is the gradient flow of the area functional and so, in principle, from a given submanifold should flow to a minimal surface. Of course, in general, the mean curvature flow develops singularities. In response ``weak solutions" of the mean curvature flow (such as the Brakke flow \cite{B}, and level set flow \cite{ES}, \cite{CGG}, and \cite{I}) have been developed. 
$\medskip$

One such approach is the mean curvature flow with surgery developed by Huisken, Sinestrari \cite{HS2} (and Brendle and Huisken \cite{BH} for the surface case) and later Haslhofer and Kleiner in \cite{HK1}). The mean curvature flow with surgery ``cuts" the manifold into pieces with very well understood geometry and topology and for this and the explicit nature of the flow with surgery is particularly easy to understand (and makes it a useful tool to understand the topology of the space of applicable hypersurfaces; see \cite{BHH} or \cite{Mra}). To be able to do this however unfortunately boils down eventually to understanding the nature of the singularities very well and establishing certain quite strong estimates, and all this has only been carried out (in $\R^{n+1}$ at least) for 2-convex compact hypersurfaces. 
$\medskip$

However the necessary estimates of Haslhofer and Kleiner as they mentioned in their paper are local in nature so can be expected to be localized in some cases. In this paper we study when this is possible, using the pseudolocality estimates of Chen and Yin \cite{CY}, and use the localized mean curvature flow with surgery to understand the level set flow - the localized mean curvature flow with surgery converges to the level set flow in a precise sense as the surgery parameters degenerate in correspondence with the compact 2-convex case. Using this we show a regularity for the level set flow and a stability result for the plane under the level set flow, showing that the mean curvature flow can be fruitfully used to study the level set flow that as far as the author knows were previously unknown. The first theorem we show in this article is a general short time existence theorem for a localized flow with surgery: 
\begin{thm}\label{exist} $\textbf{(Short time existence of flow with localized surgery)}$ Suppose $M$ is $\alpha$ noncollapsed and  $\beta$ 2-convex in an open neighborhood $U_{\Omega}$ of a bounded open set $\Omega$, and that there is $\delta, C > 0$ for which it can be guaranteed $|A|^2 < C$ in the complement of $\Omega_t$ on the time interval $[0,\delta]$ for any piecewise smooth mean curvature flow starting from $M$, with the discontinuities only occuring within $\Omega_t$. Then there exists $\eta > 0$, $\eta \leq \delta$, so that $M$ has a flow with surgery and is $\hat{\alpha} < \alpha$ non collapsed and $\hat{\beta} < \beta$ 2-convex within $\Omega_t$ on $[0, \eta]$. 
\end{thm} 
The stipulation concerning the singularities at first glance might seem rather restrictive perhaps but actually this can be guaranteed by pseudolocality estimates that control the curvature of a point through a flow just by the curvature at nearby points - this is explained in more detail after the proof of theorem \ref{exist} in section 3. 
$\medskip$

It was pointed out by Lauer in \cite{L} and independently Head in \cite{Head} that, for Huisken and Sinestrari's definition of the mean curvature flow with surgery, as the surgery parameters are allowed to degenerate the corresponding flows with surgery Hausdorff converge to the level set flow as defined by Illmanen in \cite{I}. Important for the next result and as justification of the definition of localized flow with surgery we extend Lauer's methods to show:
\begin{thm} $\textbf{(Convergence to level set flow)}$ Given $M$ if there exists a mean curvature flow with surgery, as constructed in theorem \ref{exist}, on $[0,T]$, then denoting the surgery flows $(M_t)_i$ starting at $M$ where the surgery parameter $(H_{th})_i \to \infty$ as $i \to \infty$, we have the $(M_t)_i$ as sets in $\R^{n+1} \times [0,T]$ Hausdorff converge subsequentially to the level set flow $L_t$ of $M$ on $[0,T]$. 
\end{thm}
The meaning of the surgery parameter $H_{th}$ will be described in the next section. We point out here that to overcome a technical hurdle in using Lauer's method we use ideas from the recent paper of Hershkovits and White \cite{HW} - name we use a result of their's that for us gives a way to ``localize" the level set flow.  We will mainly be interested though in such hypersurfaces with surgery satisfying additional assumptions that essentially control in a precise sense how far $M$ deviates from a plane $P$:
\begin{defn} We will say $M$ is $(V,h, R, \epsilon)$ controlled above a hypersurface $P$ in a bounded region $\Omega \subset M$ when
\begin{enumerate}
\item $M \cap \Omega$ lies to one side of the hypersurface $P$ 
\item there exists $0 < V \leq \infty$ so that the measure of points bounded initially bounded by $P$ and $\Omega$ is less than $V$. 
\item  The supremum of the height of $M$ over $P$ is bounded by $0 < h  \leq \infty$. 
\item In the $R$-collar neighborhood $C_R$ of $\partial \Omega$, $M$ is graphical over $P$ with $C^4$ norm bounded by $3\epsilon$. 
\end{enumerate} 
\end{defn}

The definition above is a bit obtuse but is essentially that the flat norm of $M$ over $P$  roughly the volume discrepancy $V$ and $h$ bound how bulky $M$ is over $P$. The definition could also possibly be phrased in terms of the flat norm of $M$ over $\Omega \subset P$ but for our applications in mind we want to keep $h$ and $V$ decoupled. By $M$ lying above $P$ we mean that $M$ lies on one side of $P$ and where $M$ is graphical over $P$ its outward normal points away from $P$ (equivalently, thinking of $M$ as the boundary of a domain $K$ so that the outward normal of $M$ is pointing outside $K$, the halfspace bounded by $P$ disjoint from $M$ lies in $K$).  The statement about the $R$-collar neighborhood of $\partial \Omega$ is for an eventual use of the Brakke regularity theorem and ensures the edges of $M$ don't ``curl up" much, see below for the case $P$ is a plane (we will mainly be interested in the case the hypersurface is extremely close to a plane). 
  \begin{center}
$\includegraphics[scale = .3]{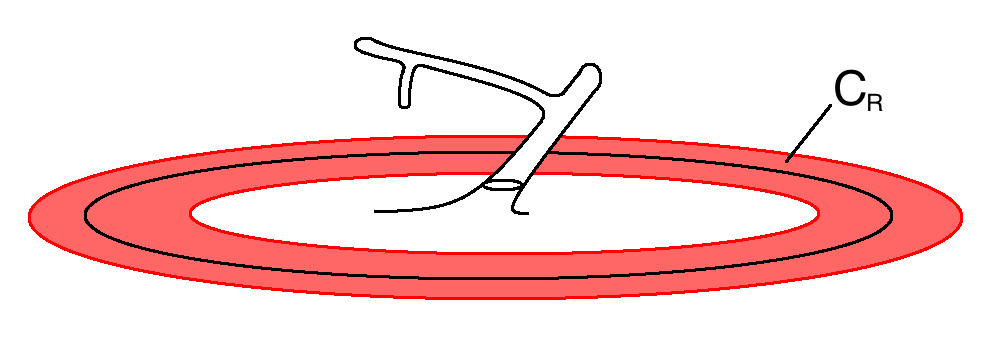}$
\end{center}

With this definition in hand, let's define the sets our regularity and stability theorems concern. The first one corresponds to the regularity result; note that for notational compactness later on we also package assumptions on $\alpha$-noncollapsedness and $\beta$ 2-convexity in $C_R$, although this could be easily modified to just concern some open set containing $\partial\Omega$: 
\begin{defn} The set $\Sigma = \Sigma( \alpha, \beta\}, \{c, S\}, \{V, h, R, \epsilon \}, \{P, \Omega\})$ is the set of hypersurfaces $M^n \subset \R^{n+1}$ satisfying:
\begin{enumerate}

\item locally $\alpha$-noncollapsed: $M \cap (\Omega \cup C_R)$ is $\alpha$-controlled in its interior. 
\item locally $\beta$-two convex: in $M \cap (\Omega \cup C_R)$ $(\lambda_1 + \lambda_2) > H\beta$.
\item supported boundary curvature: there exists $1 >> c > 0$, $S \in \{Q \in Sym(M) \mid \lambda_1(Q) + \lambda_2(Q) > 0\}$ such that $H > c\epsilon$, $A > \epsilon S$ in $C_R$
\item $M$ is $(V,h, R, \epsilon)$ controlled over the plane $P$ in the region $\Omega$.  
\end{enumerate}
\end{defn}
Note we only assume control on $\alpha$ and $\beta$ but not an initial mean curvature bound $\gamma$ (referring to the definition of an $\Balpha$-controlled domain for the surgery in Haslhofer and Kleiner's definition, see below). Items 3 and 4 imply a uniform lower bound $\eta_0$ on $\eta$ from theorem \ref{exist} for $M \in \Sigma$. With our notation and sets defined we finally state our convergence theorem; the proof crucially uses the mean curvature flow with surgery to easily get a good estimate on the height of the level set flow after a short time. 
\begin{thm}\label{smooth} $\textbf{(Local Brakke regularity type theorem for the LSF)}$ There are choices of parameters $\{\alpha, \beta\}, \{c,S\}, \{V, h, R, \epsilon \}, \{P, \Omega\}$ so that if $M \in \Sigma$ and 
\begin{enumerate}
\item has polynomial volume growth, and 
\item Is either compact or $C^0$ asymptotically flat in the sense of definition 1.3 below,
\end{enumerate}
then there is some $\eta$  on which a surgery flow of $M$ exists on $[0,\eta]$ by theorem \ref{exist}. For a given $T \in [\frac{\eta}{2}, \eta]$, there are choices of $(V,h,R, \epsilon)$ for which $L_T$ is a smooth graph over $P$.
\end{thm} 
The choice of constants $V, h$ depend on $\hat{\alpha}$, $\hat{\beta}$,  found in the existence theorem 1.1. $\hat{\alpha}, \hat{\beta},$  in turn depend on $R, \epsilon, c,$ and $S$ and also $\alpha$ and $\beta$. Since there are many parameters and their interdependence is somewhat complicated we describe explicitly after the proof of theorem \ref{smooth} in section 5 how one could choose parameters $\{ \alpha, \beta\}, \{c, S\}, \{V, h, R, \epsilon \}, \{P, \Omega\}$ so that if $M \in \Sigma$ then theorem \ref{smooth} is applicable. We will also show without much work using a general construction of Buzano, Haslhofer, and Hershkovits (theorem 4.1 in \cite{BHH}) how to construct ``nontrivial" (i.e. nongraphical, singularity forming) hypersurfaces that satisfy the assumptions of Theorem 1.3 in section 7. Such examples can also be designed to have arbitrarily large area ratios initially in a ball of fixed radius. 
$\medskip$

 Note that this theorem is an improvement on just Brakke regularity for the level set flow of $L_t$ of $M$ because we make no apriori assumptions on the densities in a parabolic ball; indeed the the hypotheses allow singularities for the LSF to develop in the regions of space-time we are considering, at which points the density will be relatively large. The point is that under correct assumptions these can be shown to ``clear out'' quickly. For a recent improvement on the Brakke regularity theorem in another, more general, direction, see the recent work of Lahiri \cite{Lah}. 
$\medskip$

This theorem is also interesting from a PDE viewpoint because the mean curvature flow is essentially a heat equation, and such result says, imagining high area ratio localized perturbations of a given hypersurface as high frequency modes of the initial condition of sorts, that in analogy to heat flow on a torus, the high frequency modes decay quickly in time. It's interesting that our arguments though use pseudolocality strongly, which is a consequence of the nonlinearity of the flow and is false for the linear heat equation. More precisely:
\begin{cor}$(\textbf{Rapid smoothing})$ Let $M$ be a smooth hypersurface with $|A|^2 < C$ for some $C > 0$ of polynomial volume growth. Suppose we perturb $M$ in some open set $U \subset M$ to get a hypersurface $\mathcal{M}$ so that:
\begin{enumerate}
\item  $\widetilde{M} \in \Sigma$, in fact that:
\item $M$ satisfies the hypotheses of theorem \ref{smooth}, and
\item $\widetilde{M}  = M$ outside set $U$
\end{enumerate}
Then then by time an appropriate $T$ as in theorem \ref{smooth} for appropriate choice of constants, $\widetilde{M}_{T}$ is smooth and has bounded curvature. 
\end{cor}
Of course taking $T$ small enough (depending on the curvature of $\widetilde{M}$ away from the perturbations) one can easily see that $\widetilde{M}_T$ is close at least in Hausdorff distance to $M$. 
$\medskip$

To state the next corollary we define a refinement of the set $\Sigma$ above, which concerns the case when $M$ is asymptotically planar with prescribed curvature decay: 
\begin{defn} The set $\Sigma_1 = \Sigma_1(\alpha, \beta\}, \{c, S\}, \{V, h, R, \epsilon\}, \{ P, \Omega\}, \{f, P_1, C_1\})$ is the set of hypersurfaces $M^n \subset \R^{n+1}$ satisfying, in addition to the set of conditions given in the definition of $\Sigma$
\begin{enumerate}
\item asymptotically planar in that in $M \cap \Omega^c$ is a graph of a function $F$ over a plane $P_1$ and furthermore writing $F$ in polar coordinates we have $||F(r, \mathbf{\theta})||_{C^2} <  f(r)$, where $f: \R_+ \to \R_+$ satisfies $\lim\limits_{r \to \infty} f(r) = 0$.
\item the hypersurface $P$ is a graph over $P_1$ with $C^1$ norm bounded by $C_1$. 
\end{enumerate}
\end{defn}

A stability statement for graphs over planes in $n \geq 3$ was noticed as a consequence of the maximum principle in an appendix of \cite{CSS} using the higher dimensional catenoids as barriers; the corollary below follows is a statement in the same spirit and follows from the flow quickly becoming graphical, the interior estimates of Ecker and Huisken, and pseudolocality:

\begin{cor}$\textbf{(Long term flow to plane)}$ With $M^n \in \Sigma_1$ asymptotically satisfying the assumptions above in theorem 1.3 in some region $\Omega$ above the origin then as $t \to \infty$ the level set $L_t$ of $M$ will never fatten and in fact will be smooth on after time $T$. It will converge smoothly to the corresponding plane $P_1$ as $t \to \infty$. \end{cor}

$\textbf{Acknowledgements:}$ The author thanks his advisor, Richard Schoen, for his advice and patience. The author also thanks the anonymous referees for their careful reading, critque, and encouragement to fill out details, which helped to much improve the clarity of the exposition. 
$\medskip$

\section{Background on the Mean Curvature Flow (With Surgery).} 
The first subsection introducing the mean curvature flow we borrow quite liberally from the author's previous paper \cite{Mra}, although a couple additional comments are made concerning the flow of noncompact hypersurfaces. The second subsection concerns the mean curvature flow with surgery as defined by Haslhofer and Kleiner in \cite{HK1} which differs from the original formulation of the flow with surgery by Huisken and Sinestrari in \cite{HK1} (see also \cite{BH}). Namely the discussion there of surgery differs from the corresponding section in the author's previously mentioned article. 
\subsection{Classical formulation of the mean curvature flow}
In this subsection we start with the differential geometric, or ``classical," definition of mean curvature flow for smooth embedded hypersurfaces of $\R^{n+1}$; for a nice introduction, see \cite{Mant}. Let $M$ be an $n$ dimensional manifold and let $F: M \to \R^{n+1}$ be an embedding of $M$ realizing it as a smooth closed hypersurface of Euclidean space - which by abuse of notation we also refer to $M$. Then the mean curvature flow of $M$ is given by $\hat{F}: M \times [0,T) \to \R^{n+1}$ satisfying (where $\nu$ is outward pointing normal and $H$ is the mean curvature):
\begin{equation}
\frac{d\hat{F}}{dt} = -H \nu, \text{ } \hat{F}(M, 0) = F(M)
\end{equation} 
(It follows from the Jordan separation theorem that closed embedded hypersurfaces are oriented). Denote $\hat{F}(\cdot, t) = \hat{F}_t$, and further denote by $\mathcal{M}_t$ the image of $\hat{F}_t$ (so $M_0 = M$). It turns out that (2.1) is a degenerate parabolic system of equations so take some work to show short term existence (to see its degenerate, any tangential perturbation of $F$ is a mean curvature flow). More specifically, where $g$ is the induced metric on $\mathcal{M}$:
\begin{equation}
 \Delta_g F = g^{ij}(\frac{\partial^2 F}{\partial x^i \partial x^j} - \Gamma_{ij}^k \frac{\partial F}{\partial x^k}) = g^{ij} h_{ij} \nu = H\nu
\end{equation}
There are ways to work around this degeneracy (as discussed in \cite{Mant}) so that we have short term existence for compact manifolds.
$\medskip$

For noncompact hypersurfaces $M$ in $\R^N$ with uniformly bounded second fundamental form (i.e. there is some $C > 0$ so that $|A|^2 < C$ at every point on $M$), one may solve the mean curvature flow within $B(0,R) \cap N$; by the uniform curvature bound there is $\epsilon > 0$ so that $N \cap B(0,R)$ has a mean curvature flow on $[0, \epsilon]$. Then one may take a sequence $R_i \to \infty$ and employ a diagonalization argument to obtain a mean curvature flow for $M$; the flow of $M$ we constructed is in fact unique by Chen and Yin in \cite{CY} (we will in fact use estimates from that same paper below). Since all noncompact hypersurfaces of interest will have asymptotically bounded geometry, we will always have a mean curvature flow of them for at least a short time. 
$\medskip$

Now that we have established existence of the flow in cases important to us, let's record associated evolution equations for some of the usual geometric quantities: 
\begin{itemize}
\item $\frac{\partial}{\partial t} g_{ij} = - 2H h_{ij}$
$\medskip$

\item $\frac{\partial}{\partial t} d\mu = -H^2 d\mu$
$\medskip$

\item $\frac{\partial}{\partial t} h^i_j = \Delta h^i_j + |A|^2 h^i_j$
$\medskip$

\item $\frac{\partial}{\partial t} H = \Delta H + |A|^2 H$
$\medskip$

\item $\frac{\partial}{\partial t} |A|^2 = \Delta |A|^2 - 2|\nabla A|^2  + 2|A|^4$
\end{itemize}
So, for example, from the heat equation for $H$ one sees by the maximum principle that if $H > 0$ initially it remains so under the flow. There is also a more complicated tensor maximum principle by Hamilton originally developed for the Ricci flow (see \cite{Ham1}) that says essentially that if $M$ is a compact manifold one has the following evolution equation for a tensor $S$: 
\begin{equation}
\frac{\partial S}{\partial t} = \Delta S + \Phi(S)
\end{equation} 
and if $S$ belongs to a convex cone of tensors, then if solutions to the system of ODE
\begin{equation}
\frac{\partial S}{\partial t} = \Phi(S)
\end{equation}
stay in that cone then solutions to the PDE (2.2) stay in the cone too (essentially this is because $\Delta$ ``averages"). So, for example, one can see then that convex surfaces stay convex under the flow very easily this way using the evolution equation above for the Weingarten operator. Similarly one can see that $\textbf{2-convex hypersurface}$ (i.e. for the two smallest principal curvatures $\lambda_1, \lambda_2$, $\lambda_1 + \lambda_2 > 0$ everywhere) remain 2-convex under the flow.
$\medskip$

Another important curvature condition in this paper is $\alpha \textbf{ non-collapsing}$: a mean convex hypersurface $M$ is said to be 2-sided $\alpha$ non-collapsed for some $\alpha > 0$ if at every point $p\in M$, there is an interior and exterior ball of radius $\alpha/H(p)$ touching $M$ precisely at $p$. This condition is used in the formulation of the finiteness theorem. It was shown by Ben Andrews in \cite{BA} to be preserved under the flow for compact surfaces. (a sharp version of this statement, first shown by Brendle in \cite{Binsc} and later Haslhofer and Kleiner in \cite{HK3}, is important in \cite{BH} where MCF+surgery to $n=2$ was first accomplished). Very recently it was also claimed to be true for non-compact hypersurfaces by Cheng in \cite{Ch}.
$\medskip$

Finally, perhaps the most geometric manifestation of the maximum principle is that if two compact hypersurfaces are disjoint initially they remain so under the flow. So, by putting a large hypersphere around $\mathcal{M}$ and noting under the mean curvature flow that such a sphere collapses to a point in finite time, the flow of $\mathcal{M}$ must not be defined past a certain time either in that as $t \to T$, $\mathcal{M}_t$ converge to a set that isn't a manifold.  Note this implies as $t \to T$ that $|A|^2 \to \infty$ at a sequence of points on $\mathcal{M}_t$; if not then we could use curvature bounds to attain a smooth limit $\mathcal{M}_T$ which we can then flow further, contradicting our choice of $T$. Of course this particular argument doesn't work in the noncompact case but it is easy to see using the Angenent's torus \cite{Ang} as a barrier that singularities can occur along the flow of noncompact hypersurfaces as well:
  \begin{center}
$\includegraphics[scale = .6]{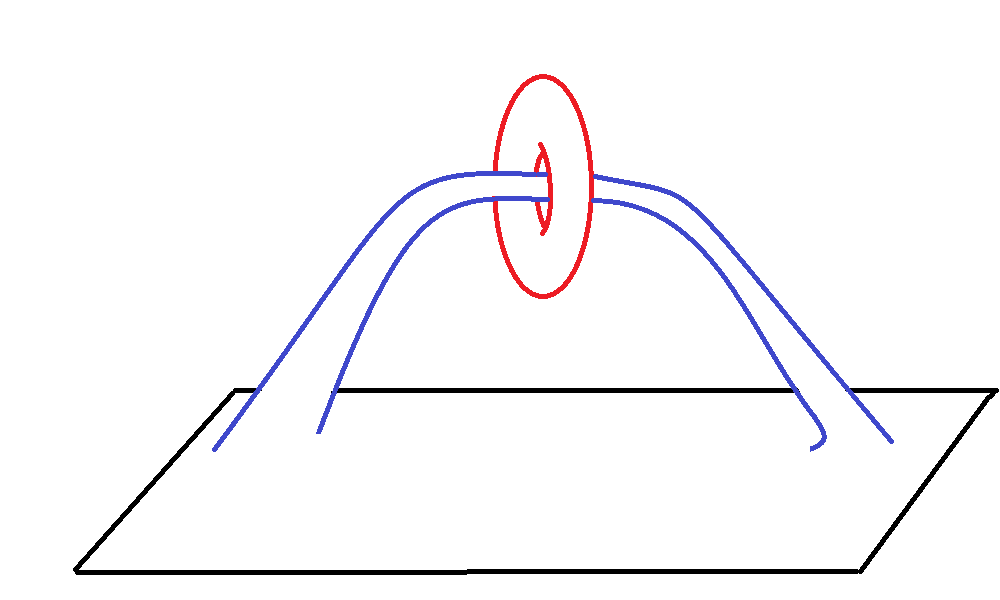}$
\end{center}

Thus as in the compact case to use mean curvature flow to study noncompact hypersurfaces one is faced with finding a way to extend the flow through singularities. Thus weak solutions to the flow are necessitated. One such weak solution is the Brakke flow, developed in Brakke's thesis \cite{B}, where a weak solution to the flow is defined in terms of varifolds. For this paper it suffices to say that the classical MCF and the LSF (defined below) are Brakke flows, and that if the density ratios of a Brakke flow are sufficiently close to 1 in a parabolic cylinder then the varifolds are actually smooth with bounded curvature within a certain time interval (this is Brakke's regularity theorem - we will be more precise about its statement in section 3) - in addition to Brakke's thesis see \cite{W1} or \cite{Lah1}. Another type of weak solution which came later is mean curvature flow with surgery: 

$\medskip$

\subsection{Mean curvature flow with surgery for compact 2-convex hypersurfaces in $\R^{n+1}$}
First we give the definition of $\Balpha$ controlled:
\begin{defn}(Definition 1.15 in \cite{HK1}) Let $\Balpha = (\alpha, \beta, \gamma) \in (0, N-2) \times (0, \frac{1}{N-2}) \times (0, \infty)$. A smooth compact domain $K_0 \subset \R^N$ is called an $\Balpha$-controlled initial condition if it satisfies the $\alpha$-noncollapsedness condition and the inequalities $\lambda_1 + \lambda_2 \geq \beta H$ and $H \leq \gamma$. 
\end{defn} 
Speaking very roughly, for the mean curvature flow with surgery approach of Haslhofer and Kleiner, like with the Huisken and Sinestrari approach there are three main constants, $H_{th} \leq H_{neck} \leq H_{trig}$. If $H_{trig}$ is reached somewhere during the mean curvature flow $M_t$ of a manifold $M$ it turns out the nearby regions will be ``neck-like" and one can cut and glue in appropriate caps (maintaining 2-convexity, etc) so that after the surgery the result has mean curvature bounded by $H_{neck}$. The high curvature regions have well understood geometry and are discarded and the mean curvature flow with surgery proceeds starting from the low curvature leftovers. Before stating a more precise statement we are forced to introduce a couple more definitions. First an abbreviated definition of the most general type of piecewise smooth flow we will consider: 
\begin{defn} (see Definition 1.3 in \cite{HK1}) An $(\alpha, \delta)-\textit{flow}$ $\mathcal{K}$ is a collection of finitely many smooth $\alpha$-noncollapsed flows $\{K_t^i \cap  U\}_{t \in [t_{i-1}, t_i]}$, $(i = 1, \ldots k$; $t_0 < \ldots t_k$) in a open set $U \subset \R^N$, such that:
\begin{enumerate}
\item for each $i = 1, \ldots, k-1$, the final time slices of some collection of disjoint strong $\delta$-necks (see below) are replaced by standard caps, giving a domain $K^{\#}_{t_i} \subset K^i_{t_i} =: K_{t_i}^-$. 
$\medskip$

\item the initial time slice of the next flow, $K_{t_i}^{i+1} =: K_{t_i}^{+}$, is obtained from $K_{t_i}^\#$ by discarding some connected components.
\end{enumerate} 
\end{defn}
For the definition of standard caps and the cutting and pasting see definitions 2.2 and 2.4 in \cite{HK1}; their name speaks for itself and the only important thing to note is that cutting and then pasting them in will preserve the $\Balpha$-control parameters on the flow. We will however give the definition of $\delta$-strong neck; below $s$ is a scaling parameter that need not concern us:
\begin{defn} (Definition 2.3 in \cite{HK1}) We say than an $(\alpha, \delta)$-flow $\mathcal{K} = \{K_t \subset U\}_{t \in I}$ has a strong $\delta$-neck with center $p$ and radius $s$ at time $t_0 \in I$, if $\{s^{-1} \cdot (K_{t_0 + s^2 t} - p) \}_{t \in (-1, 0]}$ is $\delta$-close in $C^{[1/\delta]}$ in $B_{1/\delta}^U \times (-1, 0]$ to the evolution of a solid round cylinder $\overline{D}^{N-1} \times \R$ with radius 1 at $t = 0$, where $B^U_{1/\delta} = s^{-1} \cdot ((B(p, s/\delta) \cap U) - p) \subset B(0, 1/\delta) \subset \R^N$. 
\end{defn}  

We finally state the main existence result of Haslhofer and Kleiner; see theorem 1.21 in \cite{HK1}
\begin{thm} (Existence of mean curvature flow with surgery). There are constants $\overline{\delta} = \overline{\delta}(\Balpha) > 0$ and $\Theta(\delta) = \Theta(\Balpha, \delta) < \infty$ ($\delta \leq \overline{\delta}$) with the following significance. If $\delta \leq \overline{\delta}$ and $\mathbb{H} = (H_{trig}, H_{neck}, H_{th})$ are positive numbers with $H_{trig}/H_{neck}, H_{neck}/h_{th}, H_{neck} \geq \Theta(\delta)$, then there exists an $(\Balpha, \delta, \mathbb{H})$-flow $\{K_t\}_{t \in [0, \infty)}$ for every $\Balpha$-controlled initial condition $K_0$. 
\end{thm} 
The most important difference for us (as will be evident below) between Huisken and Sinestrari's approach and Haslhofer and Kleiner's approach is that Huisken and Sinestari estimates are global in nature whereas Haslhofer and Kleiner's estimates are local. Namely, if within a spacetime neighborhood $U \times [0,T]$ it is known that the flow is uniformly $\alpha$-noncollapsed and $\beta$ 2-convex with bounded initial curvature, there are parameters $H_{th} < H_{neck} < H_{trig}$ for which surgeries can be done at times when $H = H_{neck}$ and so the postsurgery domain has curvature comparable to $H_{neck}$.

\section{Localizing the Mean Curvature Flow with Surgery.}

Recall surgery is defined for two-convex compact hypersurfaces in $\R^{n+1}$.  However, many of Haslhofer and Kleiner's estimates are local in nature and their mean curvature flow with surgery can be localized, as long as the high curvature regions (where $H > \frac{1}{2}H_{th}$ say) where surgery occurs are uniformly 2-convex (for a fixed choice of parameters). 

The main technical point then to check in performing a ``localized" mean curvature flow with surgery is ensuring that the regions where we want to perform surgeries are and remain for some time uniformly 2-convex in suitably large neighborhoods of where singularities occur. The key technical result to do so (at least for this approach) is the pseudolocality of the mean curvature flow.
$\medskip$

Pseudolocality essentially says that the mean curvature flow at a point, at least ``short term" is essentially controlled by a neighborhood around that point and that points far away are essentially inconsequential - this is in contrast to the linear heat equation. It plays a crucial role in our arguments in this section (controlling the singular set) and in some arguments in the other sections. Recall the following (consequence of the) pseudolocality theorem for the mean curvature flow of Chen and Yin:
\begin{thm} (Theorem 7.5 in \cite{CY}) Let $\overline{M}$ be an $\overline{n}$-dimensional manifold satisfying $\sum\limits_{i=0}^3 | \overline{\nabla}^i \overline{R}m| \leq c_0^2$ and inj$(\overline{M}) \geq i_0 > 0$. Then there is $\epsilon >0$ with the following property. Suppose we have a smooth solution $M_t \subset \overline{M}$ to the MCF properly embedded in $B_{\overline{M}}(x_0, r_0)$ for $t \in [0, T]$ where $r_0 < i_0/2$, $0 < T \leq \epsilon^2 r_0^2$. We assume that at time zero, $x_0 \in M_0$, and the second fundamental form satisfies $|A|(x) \leq r_0^{-1}$ on $M_0 \cap B_{\overline{M}}(x_0, r_0)$ and assume $M_0$ is graphic in the ball $B_{\overline{M}}(x_0, r_0)$. Then we have 
\begin{equation}
|A|(x,t) \leq (\epsilon r_0)^{-1}
\end{equation}
for any $x \in B_{\overline{M}}(x_0, \epsilon r_0) \cap M_t$, $t \in [0,T]$. 
\end{thm}
Of course, when the ambient space is $\R^N$, we may take $i_0 = \infty$ and since it is flat we may take $c_0 = 0$. We also highlight the following easy consequence of pseudolocality.
\begin{rem} If there are in addition initial bounds for $|\nabla A|$ and $|\nabla^2 A|$ then we also obtain bounds on $|\nabla A|$ and $|\nabla^2 A|$ a short time in the future using Chen and Yin's theorem above in combination with applying (in small balls) lemmas 4.1 and 4.2 in \cite{BH}
\end{rem}

The first and most important place pseudolocality helps us is to keep the degeneracy of 2-convexity at bay; below the MCF is normalized so it has no tangential component: 

\begin{prop} Suppose that $\Omega \subset M$ is a region in which $M$ is $\alpha$-noncollapsed  and $H > \epsilon$ on $\Omega$. Then there is $T > 0$, $\hat{\alpha} > 0$, and $\hat{\beta} > 0$ so that $\Omega$ is $\hat{\alpha}$ is non collapsed and $\hat{\beta}$ 2-convex on $[0,T]$ or up to the first singular time $T_{sing}$, if $T_{sing} <T$. 
\end{prop}
Before starting we remark that it will be clear from the proof that $T$ depends only on $|\nabla^i A|^2$, $i$ from 0 to 2 (this is coming from using the pseudolocality theorems) and a lower bound on $H$ and a lower bound on $A$ (as a symmetric matrix) in a neighborhood of $\partial \Omega$
$\medskip$

\begin{pf}

Recall the evolution equation for $H$ under the flow, that $\frac{dH}{dt} = \Delta H + |A|^2 H$. One sees by the maximum principle if $H(x)$ is a local minimum then $\frac{dH}{dt}(x) \geq 0$. This tells us that regions where $H < 0$ can't spontaneously form within mean convex regions, and in addition that for any $c$ if $\inf\limits_{x \in \Omega} H(x) > c$ in $\Omega$ intially and $H(x) > c$ on $\partial \Omega$ on $[0, T]$ then $H > c$ on all of $\Omega$ on $[0,T]$. 
$\medskip$

Let us say that $x, y \in M$ are $\alpha$-noncollapsed with respect to each other if $H(x), H(y) > 0$ and $y \not\in B(x + \nu\frac{\alpha}{H(x)}, \frac{\alpha}{H(x)})$ and vice versa. We recall from Andrew's proof \cite{BA} that provided $H > 0$, $x$ and $y$ in $M_t$ are $\alpha$-noncollapsed with respect to each other if the following quantity\footnote{Different from Andrews, we decorated our notation with $\alpha$ since this value is subject to change in our argument} is positive:
 \begin{equation}
 Z_{\alpha}(x,y,t) = \frac{H(x, t)}{2} || X(y,t) - X(x,t)||^2 + \alpha \langle X(y,t) - X(x,t), \nu(x,t) \rangle
 \end{equation} 
Of course $M_t$ is $\alpha$-noncollapsed in $\Omega_t$ if every pair of points in $\Omega_t$ is $\alpha$-noncollapsed with respect to each other. Andrews showed for closed mean convex hypersurfaces that $\alpha$-noncollapsing was preserved by the maximum principle. He calculated that on a smooth compact manifold with respect to special coordinates about extremal points $x$ and $y$ (see above lemma 4 in \cite{BA}) the following holds.  
\begin{equation}
\begin{split}
\frac{\partial Z_\alpha}{\partial t} = \sum\limits_{i,j = 1}^n \left(g^{ij}_x \frac{\partial^2Z}{\partial x^i \partial x^j} + g^{ij}_y  \frac{\partial^2Z}{\partial y^i \partial y^j}  + 2 g_x^{ij} g_y^{jl} \langle \partial_k^x, \partial_l^y \rangle  \frac{\partial^2Z}{\partial x^i \partial y^j}\right) \\ + \left(|h^x|^2 + \frac{4 H_x(H_x - \alpha h_{nn}^x)}{\alpha^2} \langle w, \partial_n^y \rangle^2 \right) Z
\end{split}
\end{equation} 
We see then that checking at values of $x$ and $y$ which minimize $Z_\alpha$ the second derivative terms are positve, so that if $Z_\alpha$ is initially nonnegative it stay so. In the our case we are interested in the noncollapsedness of a set with boundary, $\Omega_t$, but we see if for a time interval $[0,T]$ we can show there is an $\hat\alpha \leq \alpha$ so tht if $x$ and $y$ are points that minimize $Z_{\hat\alpha}$ they must be within the interior of $\Omega$, then the same argument will go through to show $\hat\alpha$ noncollapsing is preserved under the flow (note that if a set is $\alpha$-noncollapsed and $\hat\alpha \leq \alpha$, it is also $\hat\alpha$ noncollapsed). 
$\medskip$

To do this, note that in the definition 1.3 for our set $\Omega$ we have $\alpha$-noncollapsing in a neighborhood $U_\Sigma$ of $\Omega$. By pseudolocality $|A|^2$ at every point $p \in \partial \Omega$ will be bounded, for a short time, by some constant just depending on initial bounds of $|A|^2$ in a neighborhood. Since $M$ is initially smooth and $\Omega$ is bounded there are apriori uniform bounds on $|\nabla A|$, $|\nabla^2 A|$ in a neighborhood of $\partial \Omega$, which remark 3.1 above implies in combination with the last sentence impies there are uniform bounds on these quantities a short time later along $\partial \Omega$ just depending on the initial data. 
$\medskip$

Since the evolution equation $\frac{dH}{dt} = \Delta H + |A|^2 H$ is bounded by combinations of $|A|, |\nabla A|$, and $|\nabla^2 A|$ there is thus on some small forward time a uniform bound on $\frac{dH}{dt}$. Thus there is a $T >0$ just depending on $|A|, |\nabla A|$, $|\nabla^2 A|$ and $c$ for which $H(p) > c/2$ on $[0,T]$ for $p \in \partial \Omega$.  Also as a consequence of pseudolocality we see in a suitable half collar neighborhood $V$ of $\partial \Omega$ interior to $\Omega$ (see the figure below) the curvature is bound on $[0,T_1]$ by say $C$ (potentially huge). 
  \begin{center}
$\includegraphics[scale = .7]{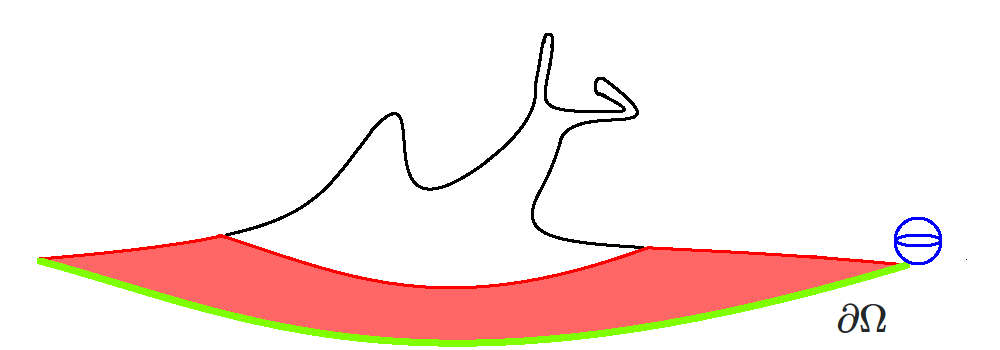}$
\end{center}

We see then there is an $\overline{\alpha} \leq \alpha$, for which we may ensure that spheres osculating $\partial \Omega_t$ of radius $\frac{\overline{\alpha}}{\epsilon}$ don't touch points in $\Omega_t$ on $[0,T_1]$. Taking $\hat\alpha = \overline\alpha/2$, we see as discussed above that $\Omega_t$ must be $\hat\alpha$ noncollapsed on $[0, T_1]$. 
$\medskip$

$\beta$-noncollapsedness is a pointwise inequality and that there is such a $\hat{\beta}$ on some fixed time $[0, T_2]$, $T_2 \leq T_1$, follows by pseudolocality as with mean convexity explained above.  
\end{pf}

Note the above theorem had no stipulation on the curvature far in the interior of $\Omega$. We are now ready to prove the short term existence theorem for the flow with surgery:
$\medskip$

\begin{pf}(of theorem \ref{exist}) Proposition 3.2 above yelds a time interval $[0, T]$ and constants $\hat\alpha, \hat\beta$ in which the set $\Omega_t$ must be $\hat\alpha$ noncollapsed and $\hat\beta$ 2-convex up to 
\begin{enumerate}
\item time $T$, if $M_t$ has a smooth flow on the interval $[0,T]$, or
\item the first singular time, which we denote $T_1$, if $T_1 < T$. 
\end{enumerate} 
In the first case, there is nothing to do. In the second case, possibly taking $T$ smaller so that $T \leq \delta$ as in the theorem statement, we know the first singularity must occur within $\Omega$ and since $T$ was choosen in the proof proposition 3.1 so that, in particular, no singularities occur along $\partial \Omega_t$, we know the singularities must be taking place in the interior of $\Omega$. By the existence of surgery for our $\hat{\alpha}$ and $\hat{\beta}$ there exists choices of parameters $H_{th} < H_{neck} < H_{trig}$ so that the surgery can be done when $H = H_{trig}$ and the curvature postsurgery will be comparable to $H_{neck}$. Furthermore the surgery parameter $H_{th}$ can be taken to be as large as one wants - we will take it larger than what $H$ could possibly obtain in $\Omega_t$ on $[0,T]$ ($H$ large implies $|A|^2$ is large). Hence a surgery can be done at a time $T^*$ before $T_1$ then, and so that $\Omega_{T^*}$ postsurgery is also $\hat\alpha$ noncollapsed and $\hat\beta$ 2-convex. 
$\medskip$

The curvature within $\Omega_{T^*}$ after the surgery will be bounded by approximately $H_{neck}$. The region outside of $\Omega_{T^*}$ will not be affected by the surgery of course, and since we stipulate we can guarantee no singularities occur outside of $\Omega_t$ on $[0, \delta]$ for any piecewise smooth flow $M_t$ starting from $M$, where the discontinuities are restricted to $\Omega_t$, the curvature on $\Omega_{t}^c$ is bounded on $[0, \epsilon]$ by some uniform constant $C'$. 
$\medskip$

Hence, we may restart the flow for some definite amount of time, if the next singular time, $T_2$, is less than $T$ we repeat the process described above. Refering to the conclusion of the theorem, $T$ will be taken to be $\eta$. 
\end{pf} 
$\medskip$

Of course, pseudolocality can be used to easily show many examples where singularities won't occur outside some fixed subset $\Omega$ for a fixed time interval of a piecewise smooth MCF that is continuous outside of $\Omega_t$. More precisely,  suppose the curvature in $U_\Omega \cup \Omega^c$ is bounded by a uniform constant, say $C_1$, and let $r_1 > 0$ be the infimum of the distance between $\partial U$ and $\partial \Omega$. Thne we see every point $x \in \Omega^c$ has a neighborhood $B(x, r_1)$ within which the curvature is bounded by $C_1$. taking $r_1$ possibly smaller, we may ensure $C_1 < \frac{1}{r_1}$. Then we can apply theorem 1.3 to see there is a time $\delta$, if all points $y \in B(x, r_1), x \in \Omega_t^c$ move by the MCF, on which the curvature at every point in $\Omega^c$ would be bounded by $\frac{1}{\epsilon r_1}$, where $\epsilon$ is the dimensional constant from theorem 3.1 - hence no singularities could occur outside of $\Omega_t$ on some short time interval. 

\section{Convergence to Level Set Flow.}
In \cite{L, Head} Lauer and Head respectively showed that as the surgery parameters degenerate, that is as $H_{th} \to \infty$, the flow with surgery Hausdorff converges to the level set flow. Strictly speaking, his theorem was for compact 2-convex hypersurfaces $M^n$, $n \geq 3$, using the surgery algorithm of Huisken and Sinestrari \cite{HS2}. As Haslhofer and Kleiner observed (see proposition 1.27 in \cite{HK1}) it is also true for their algorithm; we will show it is true for our localized surgery. 
$\medskip$

This also serves as justification for our definition of the mean curvature flow with surgery; it was important we designed our surgery algorithm to produce a weak set flow (see below). Another important observation is that, using theorem 10 of Hershkovits and White in \cite{HW}, we can ``localize" the level set flow so can get away with showing convergence to the level set flow near the singularities (in the mean convex region of $M_t$), roughly speaking.
$\medskip$

 First we record a couple definitions; these definitions are originally due to Illmanen (see \cite{I}). It is common when discussing the level set flow (so we'll do it here) to consider not $M$ but a set $K$ with $\partial K = M$ chosen so that the outward normal of $K$ agrees with that of $M$. We will quite often abuse notation by mixing $M$ and its corresponding $K$ though, the reader should be warned. When $M$ is smooth the flow of $K$ is just given by redefining the boundary of $K$ by the flow of $M$. 

\begin{defn} (Weak Set Flow). Let $W$ be an open subset of a Riemannian manifold and consider $K \subset W$. Let $\{\ell_t\}_{t \geq 0}$ be a one -parameter family of closed sets with initial condition $\ell_0 = K$ such that the space-time track $\cup(\ell_t \times \{t\}) \subset W$ is relatively closed in $W$. Then $\{\ell_t\}_{t \geq 0}$ is a weak set flow for $K$ if for every smooth closed surface $\Sigma \subset W$ disjoint from $K$ with smooth MCF defined on $[a,b]$ we have 
\begin{equation}
\ell_a \cap \Sigma_a = \emptyset \implies \ell_t \cap \Sigma_t = \emptyset
\end{equation} 
for each $t \in [a,b]$
\end{defn} 
The level set flow is the maximal such flow:
\begin{defn} (Level set flow). The level set flow of a set $K \subset W$, which we denote $L_t(K)$, is the maximal weak set flow. That is, a one-parameter family of closed sets $L_t$ with $L_0 = K$ such that if a weak set flow $\ell_t$ satisfies $\ell_0 = K$ then $\ell_t \subset L_t$ for each $t \geq 0$. 
The existence of a maximal weak set flow is verified by taking the closure of the union of all weak set flows with a given initial data. If $\ell_t$ is the weak set flow of $K \subset W$, we denote by $\hat{\ell}$ the spacetime track swept out by $\ell_t$. That is
\begin{equation}
\hat{\ell} = \bigcup\limits_{t \geq0} \ell_t \times \{t\} \subset W \times \R_+
\end{equation}
\end{defn}

\begin{rem} Evans-Spruck and Chen-Giga-Goto defined the level set flow as viscosity solutions to 
 \begin{equation} 
 w_t = | \nabla w | \text{Div} \left( \frac{\nabla w}{|\nabla w|} \right)
 \end{equation}
 but one can check (see section 10.3 in \cite{I}) that this is equivalent to the definition we gave above. 
\end{rem}

Theorem 1.2, stated more precisely then:
\begin{thm} (convergence to level set flow) Let $M \subset \R^{n+1}, n \geq 2$ be so $M$ has mean curvature flows with surgery $(M_t)_i$ as defined above on $[0,T]$ where $(H_{th})_i \to \infty$. Then
\begin{equation}
\lim\limits_{i \to \infty} \hat{(M_t)_i} = \hat{L_t}
\end{equation} 
in Hasudorff topology. 
\end{thm}

The argument of Lauer strongly uses the global mean convexity of the surfaces he has in question; in our case we only have two convexity in a neighborhood about the origin though. To deal with this we recall the following theorem of Hershkovits and White we had mentioned before:
\begin{thm} (Theorem 10 in \cite{HW}) Suppose that $Y$ and $Z$ are bounded open subsets of $\R^{n+1}$. Suppose that $t \in [0,T] \to M_t$ is a weak set flow of compact sets in $Y \cup Z$. Suppose that there is a continuous function 
\begin{center}
$w: \overline{Y \cup Z} \to \R$
\end{center}
with the following properties:
\begin{enumerate}
\item $w(x,t) = 0$ if and only if $x \in M_t$
\item For each $c$,
\begin{center}
$t \in [0,T] \to \{x \in Y : w(x,t) = c \}$
\end{center}
defines a weak set flow in $Y$. 
\item $w$ is smooth with non-vanishing gradient on $\overline{Z}$
\end{enumerate}
Then $t \in [0,T] \to M_t$ is the level set flow of $M$ in $\R^{n+1}$
\end{thm} 
Before moving on, a remark on applying the theorem to above to all the situations encountered in this article:
\begin{rem} Its clear from the proof of the theorem above that the theorem will also hold if the level sets have bounded geometry away from the surgery regions (so as to obtain the bounds in the paragraph above equation (12) in \cite{HW}). In this case $Z$ need not be bounded. In particular, the result above holds for asymptotically planar hypersurfaces - a more general class of hypersurfaces for which this is true certainly seems possible as well. \end{rem}

Hershkovits and White use this theorem to show that flows with only singularities with mean convex neighborhoods are nonfattening - previously this was only known for  hypersurfaces satisfying some condition globally like mean convexity or star shapedness. They use the theorem above to ``localize" the level set flow by interpolating between two functions of nonvanishing gradient; the distance function to the mean curvature flow of $M$ near the smooth regions and the arrival time function near the singular set (the mean convexity ensures the arrival time function has nonvanishing gradient). For our case it essentially means we only need to prove convergence of the level set flow in the mean convex region, where the singularities are stipulated to form. 
$\medskip$

So let's prove the local convergence in the mean convex region - we proceed directly as in \cite{L}. First we note the following (see lemma 2.2 in \cite{L}). Denote by $(M_{H})_t$ to be the mean curvature flow with surgery of $M$ with surgery parameter $H_{trig} = H$:
\begin{lem} Given $\epsilon > 0$ there exists $H_0 > 0$ such that if $H \geq H_0$, $T$ is a surgery time, and $x \in \R^{n+1}$, then 
\begin{equation}
B_\epsilon(x) \subset (M_H)_T^- \implies B_\epsilon(x) \subset (M_H)_T^+
\end{equation} 
\end{lem}
This statement follows jsut as in Haslhofer and Kleiner (again, proposition 1.27 in \cite{HK1}). To see briefly why it is true, since the necks where the surgeries are done are very thin, how thin depending on $H$, for any choice depending on $\epsilon > 0$ we can find an $H$ so that a ball of radius $\epsilon$ can't sit inside the neck. Hence any such ball must be far away from where any surgeries are happening.
$\medskip$

We see each of the $(M_H)_t$ are weak set flows since the mean curvature flow is and at surgery times $s_i$, $(M_H)_{s_i}^+ \subset (M_H)_{s_i}^-$. Hence $\lim\limits_{H \to \infty}$ is also. We see from how our surgery is defined in the bounded region $\Omega$ containing the surgeries that $M_t \cap \Omega$ is uniformly two convex on $[0,T]$, so that for $\epsilon > 0$ sufficiently small there exists $t_\epsilon > 0$ so that in $\Omega$:
\begin{equation}
d(M, M_{t_\epsilon}) = \epsilon 
\end{equation} 
Let $\Pi_\epsilon \subset \R^{n+2}$ be the level set flow of $M_{t_\epsilon}$. Then $\Pi_\epsilon$ is the level set flow of $K$ shifted backwards in time by $t_\epsilon$ (ignoring $t < 0$).  Let $H_0 = H_0(\epsilon)$ be chosen as in the lemma above.
$\medskip$

$\textbf{Claim:}$ $\Pi_\epsilon \subset M_H$ in $B(0,R)$ for all $H \geq H_0$. 
Let $T_1$ be the first surgery time of $M_H$. Since $\partial K_H$ is a smooth mean curvature flow on $[0,T_1)$ and $\Pi_\epsilon$ is a weak set flow the distance between the two is nondecreasing on that interval. Thus $d((\Pi_\epsilon)_T), (\partial M_H)^-_T) \geq \epsilon$ in $\Omega$ from our choice of $\epsilon$. Applying the lemma we see this inequality holds across the surgery as well. We may then repeat the argument for subsequent surgery times.  
$\medskip$

Since $\lim\limits_{\epsilon \to 0} \Pi_\epsilon = \hat{L}$ in $B(0,R)$ the claim implies $\hat{L} \subset \lim_{H \to \infty} M_H$ in $\Omega$ since the limit of relatively closed sets is relatively closed in Hausdorff topology. On the other hand as we already noted each mean curvature flow with surgery is a weak set flow for $M$. Hence the limit is also so that $\lim_{i \to \infty} (M_i) \subset \hat{L}$ in $B(0,R)$.
$\medskip$

Away from the mean convex set by assumption we have uniform curvature bounds (in our definition of mean curavture flow with localized surgery, uniform curvature bounds are assumed to occur outside the surgery regions) so for the flows with surgery $(M_i)_t$ we can pass to a Hasudorff converging subsequence that converges smoothly away from the mean convex surgery regions, and the limit by Hershkovits and White's theorem must be the level set flow. Hence we get that globally the flows with surgery converge in Hausdorff sense to the level set flow $L_t$ of $M$.

\section{A Variant of the Local Brakke Regularity Theorem for the LSF.}
 In this section we prove theorem \ref{smooth}. For the sake of reducing notational clutter we will prove for the case $P$ is the plane $x_{n+1} = 0$ - we will then easily explain why the conclusion will also be true for convex $P$ appropriately close to a plane. Also we denote (like above) $M = \partial K$ and $L_t$ the level set flow of $M$. 
 $\medskip$
 
By (4) in definition 1.1 we get a uniform lower bound on the time $T_a$ for which the surface does not intersect $P$ in $\Sigma$, without loss of generality in this section $\eta < T_a$. With that being said we show the following height estimate on mean convex flows with surgery in a ball:
 \begin{lem} Fix $\epsilon > 0$ and suppose $M \in \Sigma$. By theorem \ref{exist}, there is a $\eta > 0$ so that a flow with surgery exists out of $M$ - let $M_t$ be any such flow (i.e. no stipulations on $H_{trig}$). Then for any $0 <T < \eta$, there exists $V,h$ so that $M_T$ is in the slab bounded by the planes $x_{n+1} = \epsilon$ and $x_{n+1} = 0$
 \end{lem}

 \begin{rem} Note that no curvature assumptions are made so we may freely use this lemma as we let the surgery parameters degenerate. Also note since the post surgery domains (immediately after surgery) are contained in the presurgery domains, it suffices to consider smooth times for the flow. \end{rem}

\begin{pf}

 Denote by $\Phi_{\epsilon} \subset K_t$ the set of points in $K_t$ above the plane $x_{n+1} = \epsilon$.  Note since $M_t$ is mean convex $V_t$ is decreasing under the flow, and hence $V_t < V$, where $V_t$ is defined in the obvious way. Furthermore the $\alpha$ noncollapsing condition crucially relates $V$ and the mean curvature of points on $M_t \cap \Phi_\epsilon$ since at every point there is an interior osculating sphere proportional to the curvature; thus if $p \in \Pi$ and $x_{n+1}(p) > \epsilon$ there is a constant $\mu(\epsilon,c) > 0$ so that $ \frac{\mu}{H^{n+1}(p)} \leq |B(\frac{\alpha}{H})|  \leq V$ or so that $\sqrt[n+1]{\frac{\mu}{V}} < H(p)$.
   \begin{center}
$\includegraphics[scale = .7]{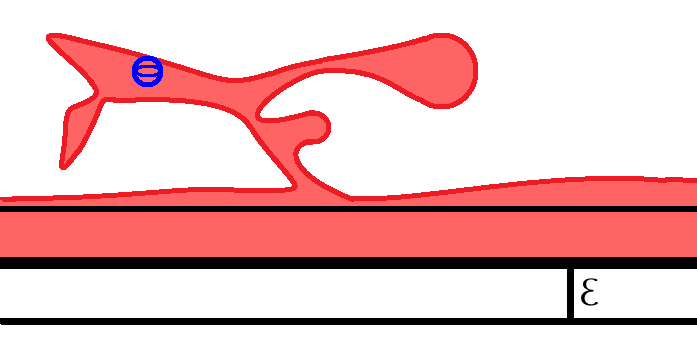}$
\end{center}
 At points on $M_t$ where the height function $x_{n+1}$ takes its maximal value the normal is pointing down, implying the height $h(t)$ of $M_t$ satisfies:
 \begin{equation}
 \frac{d h}{dt} \leq - \sqrt[n+1]{\frac{\mu}{V}}
 \end{equation} 
 we see if $h$ is small enough then the statement follows. 
 \end{pf}
 
 Note to get a simple negative lower bound for the speed of $h$ we could have also argued that there is a lower bound depend on $c$ (from definition 1.2) and $\epsilon$ as in the proof of theorem \ref{exist}- this proof is more useful though because it relates derivative of $h$ with $V$ in an explicit way. We are now ready to prove theorem \ref{smooth}; the structure of the proof is very roughly captured by the following:
 \begin{center} 
good area ratio bounds $\implies$ good Gaussian area bounds $\implies$ good Gaussian density bounds at later times $\implies$ smoothness at later times
 \end{center}
 \begin{pf} (of theorem \ref{smooth}) In theorem \ref{exist} the choice of surgery parameters didn't affect the duration of surgery (only the geometry of $\partial\Omega$ does) so we may always take a sequence of surgery flows of $[0, \eta]$ that Hausdorff converge to the level set flow of $M$ by theorem 1.2. Lemma 5.1 above also didn't depend on the choice of surgery parameters, so we thus obtain the conclusion of lemma 5.1 for the level set flow as well. With this in mind the first step is to use the lemma to get some area ratio bounds in small balls in a slab containing the plane $P$. To get these we will to use (to make our lives easier later) a slightly modified multiplicity bound theorem of White. One easily checks from the proof of the multiplicity bound theorem that it sufficed for the mean convex set $K$ to initially simply contain the slab $S$; containing the whole ball was unnecessary: 
 \begin{thm} (modified multplicity bound theorem) Let $B = B(x,r)$ be a ball, and let $S$ be a slab in $B$ of thickness $2\epsilon r$ passing through the center of the ball 
\begin{equation}
 S = \{ y \in B \mid dist(y, H) < \epsilon r \}
 \end{equation} 
 where $H$ is a hyperplane passing through the center of the ball and $\epsilon > 0$. 
Suppose $S$ is intially contained in $K$, and that $M_t \cap B$ is contained in the slab $S$. Then $K_t \cap B \setminus S$ consists of $k$ of the two connected components of $B \setminus S$, where $k$ is $0, 1$, or $2$. Furthermore
 \begin{equation} 
 area(M_t \cap B) \leq (2 - k + 2n \epsilon) \omega_n r^n. 
 \end{equation} 
 \end{thm} 
 \begin{rem} The reason it was important to not use the regular theorem (i.e. with balls) is for the sake of generality; note it is helpful for the example in section 7. 
 \end{rem}
 With this in hand we are now ready to show by time $T$ that, at least at some scales, the area ratios are very close to 1. We recall the area ratio function $\theta(L_t, x, r)$:
 \begin{equation}
 \theta(L_t, x, r) = \frac{area(L_t \cap B(x,r))}{\omega_n r^n}
 \end{equation}
Now let $M$ be as in the statement of theorem 1.3. To use the multiplicity bound theorem we first cover our plane $P$ with the ``partial'' slabs $S_p = S(p, r_0, \sigma) = \{ y \in B(p, r_0) \mid dist(y, P) < \epsilon r_0 \}$, where $r_0 > 0$ and $p \in P$. From lemma 5.1 there are appropriate choices of $V,h$ for any $0 <T_1<\eta$ so that, if initially $S(p, r_0, \epsilon) \subset  K_0$, then by time $t = T_1$, $L_{T_1}$ is contained in $\bigcup\limits_{p \in P} S_p$.
$\medskip$

Again as mentioned before the start of lemma 5.1 without loss of generality $M$ lies on one side of $P$ for $t \in [0, \eta]$. The multplicity bound theorem then holds at $t =T$ with $k = 1$, so that in each $S_x$ we have $area(L_{T_1} \cap B(x, r_0)) \leq (1 + 2n \sigma) \omega_n r_0^n$, implying that $\theta(L_{T_1}, x, r_0) < 1 + 2n \sigma$ for each $x \in P$ and in fact (since none of $L_{T_1} \cap \Omega_{T_1}$ lays outside the union of slabs) $\theta(L_{T_1}, x, r_0) < 1 + 2n \sigma$ for each $x \in \bigcup\limits_{p \in P} S_p$. Also note that in our case (i.e. nonminimal) the area ratios will not necessarily be increasing in $r$, but the control we have just at these scales is nonetheless helpful. 
 $\medskip$
 
We wish to use next the Brakke regularity theorem whch involves the Gaussian density ratio. Recall the Gaussian density ratio $\Theta(L_t, X, r)$ is given by:
 \begin{equation}
 \Theta(L_t, X, r) = \int_{y \in L_{t - r^2}} \frac{1}{(4\pi r^2)^{n/2}} e^\frac{-|y - x|^2}{4r^2} d \mathcal{H}^n y
 \end{equation} 
Where $X = (x,t)$. In analogy to area densities $\theta$ for minimal surfaces, the Gaussian density $\Theta$ are nonincreasing in $r$ along a Brakke flow, implying for a given (spatial) point $x$ that control at larger scales will give control at smaller scales forward in time. 
 $\medskip$
 
 Now we recall White's version of the local Brakke regularity theorem \cite{W1}. Since we will be considering times less than $\eta$, the flow (by Hershkovits and White, \cite{HW}) will be nonfattening so is a genuine Brakke flow. Also, as the level set flow ``biggest'' flow and nonfattening, it will agree with the Brakke flow described in section 7 of \cite{W1} and hence we will be able to apply White's version of the local Brakke regularity theorem to our case (we point this out specifically because most of that paper pertains to smooth flows up to the first singular time). 
 \begin{thm} (Brakke, White) There are numbers $\epsilon = \epsilon(N) > 0$ and $C = C(N) < \infty$ with the following property. If $\mathcal{M}$ is a Brakke flow of integral varifolds in the class $S(\lambda, m, N)$ (defined in section 7 of \cite{W1}) starting from a smooth hypersurface $M$ in an open subset $U$ of the spacetime $\R^{n+1} \times \R$ and if the Gaussian density ratios $\Theta(L_t, X, r)$ are bounded above by $1 + \epsilon$ for $0 < r< d(X,U)$, then each spacetime point $X = (x,t)$ of $\mathcal{M}$ is smooth and satisfies:
 \begin{equation}
 |A|^2 \leq \frac{C}{\delta(X, U)}
 \end{equation}
 where $\delta(X,U)$ is the infimum of $||X - Y||$ amount all spacetime points $Y \in U^c$
 \end{thm} 
 Above the statement concerning the necessary range of $r$, $0 < r< d(X,U)$, can be seen from the proof of theorem 3.1 in \cite{W1} and is a slightly stronger statement than presented in White. We also have the following important technical remark:
 
 \begin{rem} (technical remark regarding Brakke flows) Concerning the family $S(\lambda, m, N)$ of Brakke flows, it is shown in theorem 7.4 in section 7 of \cite{W1} that if $M$ is compact it has a Brakke flow in $S(\lambda, m, N)$. However one can easily check that the proof caries through if the graph of the function $u$, as in the proof, has polynomial volume growth - it is easy to construct such a function if $M$ is asymptotically flat. 
 \end{rem}
By the polynomial area growth assumption and the exponential decay of the Gaussian weight in (5.5) we see, if $R$ (in the definition of $(V,R,h,\epsilon)$-controlled) is sufficently large for a given choice of $r_1$, then $\Theta(L_{T} \cap \Omega_T^c, x, r_1)$ can be made as small as we want. By continuity of the of the Gaussian weight we see then if  for each $\delta, r_1 > 0$ we can pick $R > 0$ (so $\Omega_t$ is surrounded in a large neighborhood of nearly planar points), $r_1 > r_0 > 0$ so that if the area ratios $\theta(L_{T_1}, x, r_0) < 1 + \delta$ then $\Theta(L_{T_1}, x, r_1) < 1 +2\delta$.
$\medskip$

Now, the Brakke regularity theorem needs control over all Gaussian densities in an open set of spacetime and hence for Gaussian areas $r$ sufficiently small, but as mentioned above the favorable thing for us is that, in analogy to area densities $\theta$ for minimal surfaces, the Gaussian ratios $\Theta$ are nonincreasing in $r$ along a Brakke flow so the densities at time $T_1 + r^2$ are bounded by $1 + 4n\sigma$ for $r < r_1$. We see we had flexibility in choosing $r_0$ in the proof of the area bounds and hence we have flexibility in choosing $r_1$, so by varying $r_1$ in some small positive interval, by montonicity we obtain an interval $(a,b) \in [\frac{\eta}{2}, \eta]$ so that $\Theta(L, X, r) \leq 1 + 4n\sigma$ for all spacetime points $X \in U = \bigcup\limits_{p \in P} S_p \times (a,b)$ and $r$ sufficiently small; say $r < r^*$ for some $r^* > 0$. 
 $\medskip$
 
 Of course, the level set flow from times $(a,b)$, as these are less that $\eta$ is nonfattening (from Hershkovits and White \cite{HW}) and hence a Brakke flow, thus, taking $\sigma$ sufficiently small and $a - b$ potentially smaller (so every point $X \in U$ has $\delta(X,U) < r^*$) we may now apply Brakke regularity to get, at some time slice $T_{smooth} \in (a,b)$, uniform curvature bounds on $L_{T_{smooth}}$ a fixed distance away from $\partial \Omega$. On the other hand on the boundary we will also have curvature bounds via pseudolocality. Taking $\sigma$ smaller if need be, these curvature bounds along with the trappedness of $\Omega$ in the $\sigma$-tubuluar neighborhood of $P$ imply that $\Omega_{L_{T_{smooth}}}$ is a graph over $P$. By adjusting $T_1 \in [0,\eta]$, we may arrange $T_{smooth}$ to be any $T \in [\frac{\eta}{2}, \eta]$
 $\medskip$
 
 Now a remark regarding when $P$ is not a plane from the proofs above that if $P$ was not perfectly a plane, but merely a graph over a plane, the proof of lemma 5.1 still goes through. If the hypersurface is not a plane but so that, in the balls $B(p r_1)$ above was sufficiently close in $C_2$ norm to a plane (this depends on the $\epsilon$ necessary in the Brakke regularity theorem), we will still be able to bound the area ratios at a range of small scales as above so that the Brakke regularity theorem can be used at a later time along the flow as above. 
 \end{pf}

As promised we now discuss how one could choose parameters $\{\alpha, \beta\}, \{c, S\}, \{V, h, R, \epsilon \},$ $\{ P, \Omega \}$ so that if $M \in \Sigma$ (for these parameters) one could apply the smoothing theorem:
 \begin{enumerate}
\item Choose $\sigma$ (and hence $\epsilon$) so that the application of the Brakke regularity theorem in the above proof would hold. 
\item Having picked $\epsilon$, pick $R$ large enough so that the comment about area ratios controlling Gaussian density holds. 
\item Next, pick $\alpha, \beta, P, \Omega, V, c, S$ - these choices in particular aren't too important, although one would want $\Omega$ large compared to $R$ above to allow for topology and in section 7 the design of $P$ is important. 
\item For a given $M \in \Sigma$ theorem \ref{exist} yields a time $\eta > 0$ for which we may define a mean curvature flow with surgery.
\item Pick $h$ sufficiently small (i.e. just a bit bigger than $\epsilon$), depending on $V$, so that lemma 5.1 holds for our choice of $\epsilon$. 
\end{enumerate}
Now, one would be justifiably worried if they were concerned these sets could contain only contain trivial (i.e. already graphical) elements - it might be feared that taking $h$ small enough implies the surface is graphical for instance. However, a construction due to Buzano, Haslhofer, and Hershkovits lets us show there are nontrivial elements, as described in section 7. 
$\medskip$

Before explaining the proofs of the corollaries, we also we remark we see from the proof that $V$ and $h$ are also related and that, for a given $h > 0$, $V$ could in principle be taken sufficiently small to make the conclusion of the theorem hold - to have tall but thin 2-convex spikes. However, we see there is a lower bound on the enclosed between the plane $x_{n+1} = \epsilon$ and and $x_{n+1} = 0$ (approximately on the order $\epsilon|\Omega|$), so lemma 5.1 in practice won't be able to be used to give extremely fast speeds for $h$, refering to equation 5.1. Lemma 5.1 does work well for ``short'' spikes however. 

\section{Rapid Smoothing and LSF Long Time Convergence to a Plane Corollaries.}

\subsection{Corollary 1.4: Rapid Smoothing}
$\medskip$

This statement is essentially a ``globalization'' of Theorem \ref{smooth} and follows quickly from it. As discussed at the end of section 3, there indeed exists $\delta > 0$ on which we can ensure no singularities will occur on $[0, \delta]$ for any piecewise smooth flow starting from $M$ outside of $U$, so by theorem \ref{exist} there will be a flow with surgery localized in the open set $U$ from the statement of the corollary, on say $[0, \eta]$ where $\eta < \delta$, if the perturbations are compactly supported and $2$-convex. From theorem \ref{smooth} there is a $T < \eta$ (for appropriate choices of parameters) for which $L_t$ is smooth in $\Omega_t$ for appropriate choices of parameters. Since $T < \eta$ the surface is smooth everywhere then. As for the curvature bound, the Brakke regularity theorem gives us curvature bounds at time $T$ at fixed distance into the interior of $\Omega_T$. From definition 1.3 and the construction in theorem \ref{exist} we see the curvature will be bounded in a neighborhood of $\partial \Omega$. Putting these all together gives the statement. 
$\medskip$

\subsection{Corollary 1.5: LSF Long Time Convergence to a Plane}
Let $M \in \Sigma_1$, then theorem \ref{exist} gives some time interval $[0, \eta]$ on which the flow with surgery exists for $M$. We see by pseudolocality for a short time after the flow that, since $C_R \cup \Omega^c$ is initially a graph over the plane $P_1$ (from definition of $\Sigma_1$) this region will remain so for a short time under the flow (there will be a lower bound on this time as well for a given set of parameters). Without loss of generality then  $C_R \cup \Omega^c$ remains graphical over $P_1$ under the flow on time $[0, \eta]$. By the asymptotic planar condition the initial hypersurface $M$ (and corresponding plane $P_1$) is constrained between two parallel planes $P_1$ and $P_2$. By the avoidance principle, it must remain so under the mean curvature flow. During surgeries, high curvature pieces are discarded and caps are placed within the hull of the neck they are associated with, so $M_t$ will remain between $P_1$ and $P_2$ after surgeries as well. Thus $M_{T}$ is constrained between $P_1$ and $P_2$ for any choice of parameters and hence $L_T$ is too. 
$\medskip$

If the $M \in \Sigma_1$ for a correct choice of parameters, by time $T \in [\frac{\eta}{2}, \eta]$ $L_T \cap \Omega_T$ will be a smooth graph over $P$, and we see in the proof of theorem 1.3 by taking $\sigma$ smaller we may also arrange its Lipschtiz norm over $P$ is as small as we wish. By item (2) in definition 1.3 then for correct choices of parameters $L_T \cap \Omega_T$ will be a graph over the plane $P_1$ with bounded Lipschitz norm - hence all of $L_T$ will be by the discussion in the previous paragraph. Then we know the mean curvature flow of $L_T$, which coincides with level set flow on smooth hypersurfaces, stays graphical and its flow exists (without singularities) for all time by the classical results of Ecker and Huisken (specifically see theorem 4.6 in \cite{EH}). In fact, by proposition 4.4 in \cite{EH} one sees that as $t \to \infty$, $|A|$ and all its gradients must tend to zero uniformly in space. Since $L_T$ is bounded between two planes (this, of course, also implies its mean curvature flow is) the flow $(L_{T})_t$ of $L_{T}$ must converge to a plane. We see that in fact the plane it converges to must be $P$ since for arbitarily large times, as there will be points arbitrarily close to the plane $P$ using the asymptotically planar assumption combined with pseudolocality. 
$\medskip$

\section{Explicit Examples of Theorem 1.3}

To construct explicit nontrivial examples of mean convex regions that satisfy the hypotheses of theorem 1.3 we may use the recent gluing construction of  Buzano, Haslhofer, and Hershkovits - namely Theorem 4.1 in \cite{BHH}. It suffices to say for our purposes that it allows one to glue ``strings," tubular neighborhoods of curve segments, of arbitrarily small diameter to a mean convex hypersurface $M$ in a mean convex way. Then to construct an example, take an $\epsilon$ (from the proof of theorem \ref{smooth}) thick slab of large radius $\overline{R}$ with top and bottom parallel to the plane $x_{n+1} = 0$, which we'll denote by $S = S(\epsilon, R)$, and run it by the mean curvature flow for a very short time. The result will be convex and hence mean convex and 2-convex and remain very close to the original slab sufficiently near the origin, within say the ball $B(0,\overline{r})$. As discussed in the end of the proof of theorem \ref{smooth} the surface itself translated within this ball, translated in the $x_{n+1}$-coordinate by $- \epsilon$, can be used as the hypersurface $P$, so $h = \epsilon$ exactly with respect $P$. We see it immediately enter the associated $\epsilon$ thick slab of $P$ under the flow by mean convexity. 
$\medskip$

We then add strings to the surface near the origin; as mentioned in the proof of theorem \ref{exist} there is an $\eta > 0$ so that the surgery is possible on $[0, \eta]$ whose value doesn't depend well in the interior of $\Omega$. The strings can be taken with sufficiently small surface volume and height so that the surface must satisify lemma 5.1.
$\medskip$

By packing the strings very tightly and taking extremely small tubular neighborhoods, we can make the area ratios of $M$ in a fixed ball $B(p, \rho)$, that is the ratio of its local surface to that of the plane, as large as we want while still making the enclosed volume by the strings as small as we want (note that here outer noncollapsing isnt so important as long as nearby exterior points are intrinsically far, so a ``ball of yarn'' works). By adding small beads along the strings (that is, applying the gluing construction to glue tiny spheres along the strings in a 2-convex way) one can see using a barrier argument with the Angenent torus there are many examples of surfaces in these classes that develop singularities as well.

\end{document}